\documentclass[12pt]{article}

\usepackage{amsmath}
\usepackage{amssymb}
\usepackage{amsfonts}
\usepackage{amsthm}

\newcommand{\R}{\mathbb{R}}
\newcommand{\CC}{\mathbb{C}}

\newcommand{\N}{\mathbb{N}}

\newcommand{\Kt}{\mathfrak{H}}

\newcommand{\Id}{\mathop{\mathit{I}}}
\newcommand{\cD}{\mathcal{D}}

\newtheorem{theorem}{Theorem}
\newtheorem{Lemma}{Lemma}

\newtheorem{Remark}{Remark}

\begin{document}

\title{Multivariate Bessel functions and multivariate Hankel transforms}


\author{Victor G. Zakharov\thanks{
Institute of Continuous Media Mechanics UB RAS,
Perm 614013, Russia \&
Perm National Research Polytechnic University,
Perm 614990, Russia.
\tt victor@icmm.ru}}

\maketitle{}

\paragraph{Abstract}

The generalization, similarly to exponential multivariate bases in the Fourier transform, 
of the Bessel functions 
to
many dimensions
is offered. 
Analogously to the Fourier transform property under the differentiation, the similar
Hankel 
transform property is extended to many dimensions.

\medskip





\section{Introduction}

It is obvious fact that the 
Helmholtz differential equation
\begin{equation*}
  \frac{d^2e^{ix}}{dx^2}+e^{ix}=0,\quad x\in\R,
\end{equation*}
can be extended to any dimension $n\in\N$:
\begin{equation*}
  \sum_{j=1}^n\frac{\partial^2e^{ix}}{\partial x_j^2}+n e^{ix}=0,\quad
      x=(x_1,\dots,x_n)\in\R^n.
\end{equation*}


On the other hand,
for any {\em Bessel
function} $J_\alpha(x)$ 
\cite{123},  which is a solution of
the corresponding {\em Bessel differential equation}
\begin{equation}
  \label{BesselDifEq}
  x^2\frac{d^2u}{dx^2}+x\frac{du}{dx}+\left(x^2-\alpha^2\right)u=0,
\end{equation}
there are not direct
multivariate generalizations. 
Moreover, {\em all} the 
Bessel functions are one-dimensional.

 Certainly there are many attempts (using the Jacobi-Anger expansion,  
 see~\cite{KorschklumppWitthaut}, for example) to obtain multivariate Bessel functions. 

Nevertheless, in this very short paper, we offer more direct and simple method, in fact similar to the 
multidimensional 
exponential bases construction, to build the multivariate Bessel functions.
Moreover, since our approach is similar to the Fourier; many properties of the (multidimensional) Fourier transform
can be prorogated to the Hankel transform. 



First, using the three-dimensional {\em parabolic coordinates}
$\left(\sigma,\tau,\phi\right)$
that are expressed in terms of the Cartesian coordinates as:
\begin{align*}
  x&=\sigma\tau\cos\phi,\\
  y&=\sigma\tau\sin\phi,\\
  z&=\frac12\left(\tau^2-\sigma^2\right),
\end{align*}
consider
the Laplacian
$$
  \Delta u=
  \frac{1}{\sigma^2+\tau^2}\left[\frac{1}{\sigma}\frac{\partial}{\partial\sigma}
     \left(\sigma\frac{\partial u}{\partial\sigma}\right)
     +\frac{1}{\tau}\frac{\partial}{\partial\tau}\left(\tau\frac{\partial u}{\partial\tau}\right)\right]
   + \frac{1}{\sigma^2\tau^2}\frac{\partial^2u}{\partial\phi^2}.
 $$
 It is easily shown that, by separation of variables, the corresponding Bessel differential equation is of the form
\begin{equation}
  \label{2DBE}
  \frac{1}{\sigma}\frac{\partial}{\partial\sigma}
     \left(\sigma\frac{\partial u}{\partial\sigma}\right)
     +\frac{1}{\tau}\frac{\partial}{\partial\tau}\left(\tau\frac{\partial u}{\partial\tau}\right)
     +\left(2-\frac{\alpha^2}{\sigma^2}-\frac{\beta^2}{\tau^2}\right)u=0,\quad \alpha,\beta\in\CC;
\end{equation}
and the Bessel-related 
two-dimensional solution of equation~\eqref{2DBE} can be defined as follows
\begin{equation}
  \label{Jab}
  J_{\alpha\beta}(\sigma,\tau):=J_\alpha(\sigma)J_\beta(\tau),\qquad \sigma,\tau\in\R^+,\
  \alpha,\beta\in\CC,
\end{equation}
where $J_\alpha$, $J_\beta$ are the first kind
 Bessel functions of order $\alpha$ and $\beta$, respectively.

The second example concerns the {\em spherical Bessel functions $j_\alpha,j_\beta$} and 
the {\em spherical coordinates}.
Analogously to the previous example,  
we can state the following two-dimension differential equation
 \begin{equation}
  \label{2DSherBE}
  \frac{1}{\sigma^2}\frac{\partial}{\partial\sigma}
     \left(\sigma^2\frac{\partial u}{\partial\sigma}\right)
     +\frac{1}{\tau^2}\frac{\partial}{\partial\tau}\left(\tau^2\frac{\partial u}{\partial\tau}\right)
     +\left(2-\frac{\alpha(\alpha+1)}{\sigma^2}-\frac{\beta(\beta+1)}{\tau^2}\right)u=0,
\end{equation}
where a solution 
of equation~\eqref{2DSherBE}
is defined like~\eqref{Jab}:
\begin{equation}
  \label{SperJab}
  j_{\alpha\beta}(\sigma,\tau):=j_\alpha(\sigma)j_\beta(\tau),\qquad \sigma,\tau\in\R^+,
  \
  \alpha,\beta\in\CC.
\end{equation}

Thirdly we present the functions that are not Bessel-related.
Namely we mean the {\em associated Legendre polynomials} $P^m_l(x)$, $x\in \R$,
where the integers $l$ and $m$ are referred to as the {\em degree} and {\em order}, respectively. Hence the {\em general Legendre equation} is of the form
$$
  \frac{d}{dx}\left[\left(1-x^2\right)\frac{d}{dx}P^m_l(x)\right]+\left[l(l+1)
     -\frac{m^2}{1-x^2}\right]P^m_l(x)=0.
$$
As before, introduce the two-dimensional associated Legendre polynomials
\begin{equation}
  \label{LegPol2D}
  P^{m_1m_2}_{l_1l_2}(x,y):=P^{m_1}_{l_1}(x)P^{m_2}_{l_2}(y)
\end{equation}
that are solutions to the two-dimensional associated Legendre equation
\begin{equation}
  \label{LeqEq2D}
\begin{aligned}
  \frac{d}{dx}&\left(\left(1-x^2\right)
    \frac{dP^{m_1m_2}_{l_1l_2}}{dx}\right)
     +\frac{d}{dy}\left(\left(1-y^2\right)\frac{dP^{m_1m_2}_{l_1l_2}}{dy}\right)\\
  &+\left[l_1(l_1+1)+l_2(l_2+1)-\frac{m_1^2}{1-x^2}-\frac{m_2^2}{1-y^2}\right]
         P^{m_1m_2}_{l_1l_2} =0.
\end{aligned}
\end{equation}

It is not difficult to verify 
that offered functions~\eqref{Jab},
\eqref{SperJab}, and~\eqref{LegPol2D}
are 
solutions to differential equations~\eqref{2DBE}, \eqref{2DSherBE}, and~\eqref{LeqEq2D}, respectively.
The proof will be presented later.

Note that the second kind
functions $Y_\alpha$, $y_\alpha$, $Q^m_l$
can be
substituted for the Bessel functions $J_\alpha$, $j_\alpha$, and
the associated Legendre polynomials $P^m_l$ in
differential equations~\eqref{2DBE}, \eqref{2DSherBE}, and~\eqref{LeqEq2D}, respectively.

The problems demonstrated in the previous examples (see equations~\eqref{2DBE}, \eqref{2DSherBE}, and~\eqref{LeqEq2D})
 can be generalized to any dimension $n\in\N$.
 In particular,
 let $x:=(x_1,\dots,x_n)$, $\alpha:=(\alpha_1,\dots,\alpha_n)$; then
the $n$-dimensional Bessel functions are defined as
\begin{equation}
\label{MMB}
  J_\alpha(x):=\prod_{j=1}^n J_{\alpha_j}(x_j);
\end{equation}
and the corresponding $n$-dimensional equations are
\begin{equation*}
  \sum_{j=1}^n \frac{1}{x_j}\frac{\partial}{\partial x_j}\left(x_j\frac{\partial
     J_\alpha(x)}{\partial x_j}\right)+\left(n-\sum_{j=1}^n\frac{\alpha_j^2}{x_j^2}\right)J_\alpha(x)=0.
\end{equation*}
For other examples, the modifications of the previous formulas are obvious.

Note also that the multidimensional functions that we offered and the presented examples
of differential equations concern with the corresponding multidimensional coordinate systems (cylindrical, spherical, and ellipsoidal),
which
are independent problems to construct.
Some solutions of the mentioned problems will be discussed elsewhere.

\section{Proof}

In this section we present a proof of our approach to construct the multidimensional
Bessel-like functions.

By definition, put
\begin{equation}
  \label{MyDefin}
  \cD_{\alpha;x}f:= \frac{1}{x}\frac{\partial}{\partial x}
     \left(x\frac{\partial f}{\partial x}\right)-\frac{\alpha^2}{x^2}f,\qquad
        x\in\R,\alpha \in\CC,
\end{equation}
where $f=f(x), x\in\R,$ is some differentiable function.

By~\eqref{MyDefin}, Bessel differential equation~\eqref{BesselDifEq}
can be rewritten as follows
\begin{equation}
  \label{FixedPoint}
  \cD_{\alpha;x}J_{\alpha}(x)=-J_{\alpha}(x),
\end{equation}
where $J_\alpha$ is the Bessel function;
and we see that $J_\alpha$ is a fixed point (up to a constant factor)
of operator~\eqref{MyDefin}. 

\begin{theorem}\label{Th}
Two-dimensional Bessel function~\eqref{Jab} is a solution of two-dimensional
Bessel differential equation~\eqref{2DBE}.
\end{theorem}


\begin{proof}
Using definition~\eqref{Jab} 
and equation~\eqref{FixedPoint},
the left-hand side of
equation~\eqref{2DBE} can be written as follows
\begin{equation*}
\begin{aligned}
  \left[\cD_{\alpha;x}+\cD_{\beta;y}\right]J_{\alpha\beta}(x,y)
  =J_\beta(y)\cD_{\alpha;x}J_\alpha(x)&+J_\alpha(x)\cD_{\beta;y}J_\beta(y)\\
  &=-2J_\alpha(x)J_\beta(y)=-2J_{\alpha\beta}(x,y).
\end{aligned}
\end{equation*}
Consequently the function $J_{\alpha\beta}(x,y)$ is a root of
equation~\eqref{2DBE}.
\end{proof}

\begin{Remark}
Of course, in the proof of Theorem~\ref{Th}, we must motivate the applicability
of the constant factor rule
$(Cf)'=Cf'$
for the statements $\cD_{\alpha;x}J_{\alpha\beta}(x,y)=J_\beta(y)\cD_{\alpha;x}J_\alpha(x)$
and $\cD_{\beta;y}J_{\alpha\beta}(x,y)=J_\alpha(x)\cD_{\beta;y}J_\beta(y)$.
So, the more strong proof of the theorem is
left to the reader.
\end{Remark}

Performing the iteration of equation~\eqref{FixedPoint}, we have
$$
  \cD_{\alpha;x}^m J_{\alpha}(x)=\left(-1\right)^m J_{\alpha}(x).
$$


\section{Multidimensional Hankel transform}

\subsection{One-dimensional Hankel Transform}

Recall that the {\em Hankel transform $\Kt_\nu$} of order $\nu$ for a function $f$ is of the form
\begin{equation}
   \label{OneDimHankTr}
  \Kt_\nu\left[f(r)\right]=\Kt_\nu\left[f(r);k\right]:=\int_0^\infty f(r)J_\nu(kr)r\,dr,
\end{equation}
where $J_\nu$ is the Bessel function of the first kind of order $\nu$.

Recall also  {a property 
of the Hankel transform:
\begin{equation}
  \label{IntProperty}
  \Kt_\nu[\cD_{\nu;r}f(r);k]=-k^2\Kt_\nu[f(r);k].
\end{equation}
where $\cD_{\nu;r}$ is an operator defined by formula~\eqref{MyDefin}.

\begin{Lemma}
\label{Omit}
Substituting $J_\nu(kr)$ for $f(r)$ in equation~\eqref{IntProperty},
we can omit the Hankel transform symbol $\Kt_\nu$ in the
left-hand and right-hand sides of the equation. Namely we have
%
\begin{equation}
  \label{NewIntProperty}
    \cD_{\nu;r}J_{\nu}(kr)=
      -k^2J_{\nu}(kr).
\end{equation}
\end{Lemma}
\begin{proof}
Substituting $J_{\nu}(kr)$ for $f(r)$ in equation~\eqref{IntProperty},
applying the Hankel transform $\Kt_\nu$ to the left-hand side of the equation, and 
%
using
the orthogonality relation for the Bessel function
\begin{equation}
  \label{BessOrth}
  \int_0^\infty x J_\alpha(ux)J_\alpha(vx)\,dx=\frac{1}{u}\delta(v-u),
\end{equation}
we obtain
\begin{equation}
  \label{ququ}
 \Kt_\nu\left[\cD_{\nu;r}J_\nu(kr);a\right]
         =-a^2\frac{1}{a}\delta(k-a).
\end{equation}
Then apply the Hankel transform $\Kt_\nu$ to the right-hand side of expression~\eqref{ququ}
\begin{equation*}
  \Kt_\nu\left[-a^2\frac{1}{a}\delta(k-a);r\right]
    =-k^2J_\nu(kr).
\end{equation*}
On the other hand, applying $\Kt_\nu$ to the left-hand side of expression~\eqref{ququ}
and using the inverse property of the Hankel transform
$\Kt_\nu\Kt_\nu=\Id$, where $\Id$ is the Identity operator,
we get
$$
 \Kt_\nu\left\{\Kt_\nu\left[\cD_{\nu;r}J_\nu(kr);a\right];r\right\} = \cD_{\nu;r}J_\nu(kr).
$$
So, equation~\eqref{NewIntProperty} is valid.
\end{proof}

\begin{Remark}
Lemma~\ref{Omit} can be proved also
for the spherical Bessel functions and (probably)
for the associated Legendre polynomials.

In particular,
the orthogonality relation
\begin{equation}
  \label{SphBessOrth}
  \int_0^\infty x^2 j_\alpha(ux)j_\alpha(vx)\,dx=\frac{\pi}{2u^2}\delta(u-v)
\end{equation}
is valid
for the spherical Bessel functions;
and
the orthogonality relation
(similar to~\eqref{BessOrth},~\eqref{SphBessOrth})
for the 
Legendre polynomials
is left to the reader.
\end{Remark}

\subsection{$n$-dimensional Hankel transform}

By definition~\eqref{MMB}, Hankel transform~\eqref{OneDimHankTr} can be generalized
to any 
dimension $n\in\N$:
\begin{equation}
  \label{NDimHankTr}
  \begin{aligned}
  &\Kt_{\nu_1,\dots,\nu_n}\left[f(r_1,\dots,r_n);
     k_1,\dots,k_n\right]\\
  &\quad :=\int_{[0,\infty)^n} f(r_1,\dots,r_{n})J_{\nu_1,\dots\nu_n}(k_1r_1,\dots,k_nr_n)
            r_1\cdots r_n\,dr_1\cdots dr_n.
  \end{aligned}
\end{equation}

Using
the ideas of the proof of Lemma~\ref{Omit}, we can
state and proof the following theorem.

\begin{theorem}\label{TheoremII}
Let $P$ be an algebraic polynomial, $J_{\alpha_1,\dots,\alpha_n}$ be the $n$-dimensional
Bessel function defined by~\eqref{MMB}, and $\cD_{\alpha_1;k_1},\cD_{\alpha_2;k_2},\dots,\cD_{\alpha_n;k_n}$
be the differential operators defined by~\eqref{MyDefin}. Then we have
\begin{equation}
  \label{IntProperty1}
\begin{aligned}
  P\left(-k_1^2,\dots,-k_n^2\right)&J_{\nu_1,\dots\nu_n}\left(k_1r_1,\dots,k_nr_n\right)\\
  &=
  P\left(\cD_{\alpha_1;r_1},\dots,\cD_{\alpha_n;r_n}\right)
    J_{\nu_1,\dots\nu_n}\left(k_1r_1,\dots,k_nr_n\right).
\end{aligned}
\end{equation}
\end{theorem}

\section{Future trends}

Here we introduce some intents in the near future.
Namely we shall try to obtain\dots
\begin{itemize}
\item
\dots the analog of Theorem~\ref{TheoremII} for odd degrees;
\item
\dots 
the matrix algorithms, based on multidimensional Bessel functions and Hankel transforms, to solve
PDE('s) with {\em polynomial} coefficients;
\item
\dots 
the multidimensional functions,
in case of other than cylindrical or spherical 
coordinate systems; 
\item
\dots the multidimensional functions, in case of other than the Laplacian differential operator, in particular, the operators of odd degree;
\end{itemize}

\end{document}